\documentclass[11pt,amscd]{amsart}
\usepackage{amsaddr}
\usepackage{graphicx,amssymb,amsmath,amsthm}

\usepackage[utf8]{inputenc}
\usepackage{epstopdf}
\usepackage{url,hyperref}	
\usepackage[boxed,linesnumbered]{algorithm2e}	
\usepackage{bm}
\usepackage[numbers]{natbib}
\usepackage{color}
\usepackage{listings}
\usepackage{array}

\usepackage[margin=1.25in]{geometry}
\newcolumntype{L}{>{$}l<{$}} 
\newcolumntype{C}{>{$}c<{$}} 

\usepackage{bm,graphicx,textcomp,bbm} 

\newcommand{\eps}{{\varepsilon}}

\renewcommand*{\Xi}{{\boldsymbol{\xi}}}

\renewcommand*{\O}{{\mathcal{O}}}
\newcommand*{\Q}{{\mathcal{Q}}}

\newcommand*{\E}{{\mathbb{E}}}
\newcommand*{\V}{{\mathbb{V}}}
\newcommand*{\R}{{\mathbb{R}}}

\newcommand*{\change}[1]{#1}

\newcommand{\pav}[2]{\left\langle #1 \right\rangle_{\varepsilon_#2}}
\newcommand{\paveps}[1]{\left\langle #1 \right\rangle_{\varepsilon}}
\newcommand{\ueps}[1]{_{\varepsilon_{#1}}}
\newcommand{\bxi}{{\xi}}
\newcommand{\q}{\mathfrak{q}}

\usepackage{amsfonts}

\title{On the validity of linear response theory in high-dimensional deterministic dynamical systems}

\author{Caroline L. Wormell \and Georg A. Gottwald}
\address{School of Mathematics and Statistics, University of Sydney, NSW 2006, Australia}
\email[C. L. Wormell and G. A. Gottwald]{georg.gottwald@sydney.edu.au {\rmfamily and} ca.wormell@gmail.com} 
\thanks{{\em Mailing address}: School of Mathematics and Statistics, University of Sydney, NSW 2006, Australia}
\thanks{C.L.W. and G.A.G. designed and performed research and wrote the paper.}

\begin{document}

\begin{abstract}
	This theoretical work considers the following conundrum: linear response theory is successfully used by scientists in numerous fields, but mathematicians have shown that typical low-dimensional dynamical systems violate the theory's assumptions. Here we provide a proof of concept for the validity of linear response theory in high-dimensional deterministic systems for large-scale observables. We introduce an exemplary model in which observables of resolved degrees of freedom are weakly coupled to a large, inhomogeneous collection of unresolved chaotic degrees of freedom. By employing statistical limit laws we give conditions under which such systems obey linear response theory even if all the degrees of freedom individually violate linear response. We corroborate our result with numerical simulations. 
\end{abstract}

\maketitle

\section{Introduction}

Linear response theory (LRT) has been a cornerstone of statistical mechanics ever since its introduction in the 1960s. When valid, it allows us to express the average of some observable when subjected to small perturbations from an unperturbed state -- the system's so called {\em{response}} -- entirely in terms of statistical information from the unperturbed system. In essence, linear response theory relies on the smoothness of the invariant measure with respect to a perturbation, in the sense that there exists a Taylor expansion of the perturbed invariant measure around the unperturbed equilibrium measure. 

The development of the theory occured in statistical mechanics in the context of thermostatted Hamiltonian systems \cite{Kubo66,Balescu1,Zwanzig,MarconiEtAl08} but found applications far beyond this realm; recent years have seen an increased interest in LRT and its applications. In particular, climate scientists have resorted to LRT to study the timely question how certain observables such as the global mean temperature or local rainfall intensities behave upon increasing the ${\rm{CO}}_2$ concentration in the atmosphere. LRT has been successfully applied to several situations with macroscopic observables in various atmospheric toy models \cite{MajdaEtAl10,LucariniSarno11,AbramovMajda07,AbramovMajda08,CooperHaynes11,CooperEtAl13}, barotropic models \cite{Bell80,GritsunDymnikov99,AbramovMajda09}, quasi-geostrophic models \cite{DymnikovGritsun01}, atmospheric models \cite{NorthEtAl93,CionniEtAl04,GritsunEtAl02,GritsunBranstator07,GritsunEtAl08,RingPlumb08,Gritsun10} and in coupled climate models \cite{LangenAlexeev05,KirkDavidoff09,FuchsEtAl14,RagoneEtAl15}.

In a separate strand of research mathematicians have tried to obtain rigorous results extending the validity of LRT to deterministic dynamical systems. There was initial success by Ruelle \cite{Ruelle97,Ruelle98,Ruelle09a,Ruelle09b} in the case of uniformly hyperbolic Axiom A systems, however the works of Baladi and colleagues undermined hopes that LRT typically holds in dynamical systems \cite{BaladiSmania08,BaladiSmania10,Baladi14,BaladiEtAl15,DeLimaSmania15}. They showed that simple dynamical systems such as the logistic map do not obey LRT but rather their invariant measure changes non-smoothly with respect to the perturbation (even considering only chaotic parameter values). This poses a conundrum: how can LRT seem to be typically valid in high-dimensional systems for macroscopic observables when structural obstacles to its validity are likely to be present in its microscopic constituents? 

To justify the validity of LRT in high-dimensional systems, scientists often invoke the {\em{chaotic hypothesis}} of Gallavotti-Cohen \cite{GallavottiCohen95a,GallavottiCohen95b} according to which a high-dimensional system behaves for all practical purposes as an Axiom A system. This invocation, however, is unjustified: even if the hypothesis is true, it does not address how the equivalent Axiom A systems of the unperturbed and the perturbed system relate to each other, which is crucial for any statement on LRT. 

In a recent paper \cite{GottwaldEtAl16} we showed that breakdown of LRT might not be detectable using uncertainty quantification when analyzing time series unless the time series is very long (exceeding 1 million data points even for simple one-dimensional systems such as the logistic map, for example) and/or the observables are sensitive to the non-smooth change of the invariant measure. Consequently, the apparent observed validity of LRT in climate science might be a finite size effect. 

Here we follow a different avenue, drawing on the fact that linear response theory can be justified \cite{Haenggi78,HairerMajda10} for stochastic dynamical systems. We argue here that certain deterministic chaotic systems have stochastic limits for macroscopic observables which implies that they are amenable to LRT. Statistical limit laws of deterministic dynamical systems have recently been proven for slow variables in multi-scale systems \cite{MelbourneStuart11,GottwaldMelbourne13c,KellyMelbourne14} and for resolved degrees of freedom in high-dimensional weakly coupled systems \cite{FordEtAl65,Zwanzig73,FordEtAl87,StuartWarren99,KupfermanEtAl02,GivonEtAl04}. In both cases the diffusive limit of the macroscopic observables relies on the central limit theorem via a summation of infinitely many weakly dependent variables. We treat here the case of weak coupling whereby distinguished resolved degrees of freedom are weakly coupled to a large {\em{heat bath}} of unresolved, dissipative microscopic degrees of freedom. The central limit theorem can be justified in this situation either for sufficiently chaotic dynamics (the case we consider here) or for a collection of randomly chosen initial conditions. We introduce here a judiciously chosen toy model which considers the worst case scenario where both the resolved and the unresolved dynamics violate LRT, when considered on their own. Our main finding is that LRT can be assured in high-dimensional systems of weak coupling type, when the macroscopic resolved variables exhibit effective stochastic dynamics and when additionally the microscopic dynamics is spatially heterogeneous.\\

The paper is organized as follows. Section~\ref{s.LRT} briefly reviews LRT.  In Section~\ref{s.model} we introduce the high-dimensional weak coupling model under consideration. Section~\ref{s.homo} considers the case when the resolved scales exhibit a diffusive limit in the thermodynamic limit of an infinite-dimensional microscopic sub-system, and we show that LRT is valid. Section~\ref{s.det} treats the case when the thermodynamic limit is deterministic and LRT is not valid for infinitely many degrees of freedom. We will see, however, that for large but finite system sizes, linear response is valid for some, albeit small, range of perturbations, and the breakdown of LRT might not be detectable in typical time series for an increasing range of perturbations. We conclude with a discussion and an outlook in Section~\ref{s.discussion}.


\section{Linear response theory}
\label{s.LRT}
Consider a family of dynamical systems $f_\varepsilon:D \to D$ on some space $D$ where the map $f_\varepsilon$ depends smoothly on the parameter $\varepsilon$ and where for each $\varepsilon$ the dynamical system admits a unique invariant physical measure $\mu_\varepsilon$. An ergodic measure is called physical if for a set of initial conditions of nonzero Lebesgue measure the temporal average of a continuous observable converges to the spatial average over this measure. LRT is concerned with the change of the average of an observable $\Psi:D\to\R$,  
\begin{align*}
\E^\varepsilon [\Psi] = \int_D \Psi\, d\mu_\varepsilon
\end{align*}
upon varying $\varepsilon$. A system exhibits {\em{linear response}} at $\eps=\eps_0$,  if the derivative
\begin{align*}
\E^{\varepsilon_0} [\Psi]^\prime := \frac{\partial}{\partial\varepsilon} \E^\varepsilon [\Psi]|_{\varepsilon_0}
\end{align*}
exists. One can then express the average of an observable of the perturbed state with $\eps=\eps_0+\delta\eps$ up to $o(\eps)$ as
\begin{align*}
\E^\varepsilon [\Psi]  \approx \E^{\varepsilon_0} [\Psi] +\delta \eps\, \E^{\varepsilon_0} [\Psi]^\prime ,
\end{align*}
which may be determined entirely in terms of the statistics of the unperturbed system and its invariant measure $\mu_{\eps_0}$ using so-called linear response formulae \cite{Ruelle09a,Ruelle98,Baladi14}. 
A sufficient condition for linear response is therefore that the invariant measure $\mu_\varepsilon$ is differentiable with respect to $\varepsilon$. If the limit does not exist, we say there is a breakdown of linear response. We assume that the observable captures sufficient dynamic information about the dynamical system; for example, an odd observable on a system symmetric about $0$ would be identically zero regardless of whether the system exhibits linear response or not.


\section{The model}
\label{s.model}
\begin{figure}
	\centering
	\includegraphics[width=0.4\linewidth,clip=true,angle=90]{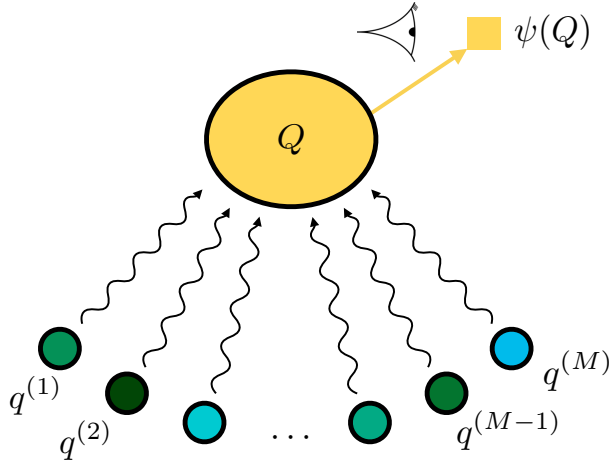}
	\caption{Sketch of the toy model set up of a macroscopic resolved variable $Q$ which is observed with observable $\psi(Q)$. $Q$ is weakly coupled to unresolved microscopic variables $q^{(j)}$, $j=1,\cdots,M$. The microscopic sub-system is heterogeneous with each microscopic variable evolving independently according to its own randomly drawn parameters, as indicated by the different coloured shadings.}
	\label{f.model}
\end{figure}

\change{
We introduce an exemplary high-dimensional toy system where each individual component does not obey linear response. We consider the case of a single resolved macroscopic degree of freedom $Q$ weakly coupled to $M$ unresolved microscopic degrees of freedom $q^{(j)}, j = 1, \ldots, M$. The microscopic dynamics is assumed to evolve independently of the macroscopic dynamics and independently of each other. Further, we make the natural assumptions that the microscopic dynamics is heterogeneous in the sense that each microscopic variable $q^{(j)}$ evolves with their own parameter $a^{(j)}$, drawn from a smooth distribution $\nu$.  To study the linear response of macroscopic observables $\Psi(Q)$ we consider perturbations of the parameters of the microscopic dynamics of the form $a^{(j)} = a_0^{(j)} +  \epsilon a_1^{(j)}$. Figure~\ref{f.model} provides a graphical illustration of our set up.}

\change{
To illustrate how a high-dimensional system can exhibit linear response for macroscopic observables even if their microscopic constituents do not obey LRT, we make the worst case assumption that neither the microscopic variables nor the macroscopic variables obey linear response when viewed in isolation. For the purposes of this paper we use the prototypical example of a logistic-type map for a dynamical system which violates LRT \cite{BaladiSmania08,BaladiSmania10,Baladi14,BaladiEtAl15,DeLimaSmania15}.}

\change{
To be specific, the macroscopic variable $Q$ evolves according to a logistic map
\begin{align}
Q_{n+1} = A\,Q_n(1-Q_n),
\label{e.Qn}
\end{align}
with parameter 
\begin{align*}
A = A_0 +A_1 Z_n
\end{align*} 
driven by the unresolved, microscopic variables through the coupling term
\begin{align}
Z_n = \frac{1}{M^\gamma}\sum_{j=1}^{M}\phi^{(j)}_n,
\label{e.Zn}
\end{align}
with scaling parameter $\gamma \geq \tfrac{1}{2}$. Here 
\begin{align*}
\phi^{(j)}_n = \phi(q^{(j)};a^{(j)})
\end{align*}
is a H\"older continuous function of the microscopic variables $q_{(j)}$. The $M$ unresolved microscopic degrees of freedom $q^{(j)}$ evolve according to modified logistic maps of the form
\begin{align}
{\small{
\left( q_{n+1}^{(j)},r_{n+1}^{(j)}\right) 
 = 
\begin{cases}
\left( q_{n}^{(j)},2r_{n}^{(j)}\right) & r_n^{(j)}<\tfrac{1}{2}\\
\left(a^{(j)}\,q_n^{(j)}(1-q_n^{(j)}),2 r_n^{(j)}-1\right) &  r_n^{(j)}\ge\tfrac{1}{2}
\end{cases},
}}
\label{e.qn}
\end{align}
each with their particular parameter $a^{(j)}$. The modification of the logistic map as a cocycle over a mixing doubling map for $r_n$ assures that the overall dynamics is mixing (thereby avoiding any periodic dynamics of the microscopic variables  $q^{(j)}$). The modified map is constructed such that its marginal invariant measure of  $q^{(j)}$ recovers exactly the physical measure of the standard logistic map with the same parameter $a^{(j)}$. Hence the microscopic dynamics (\ref{e.qn}) violates LRT while being chaotic.} 

\begin{figure}
	\centering
	\includegraphics[width=0.7\linewidth]{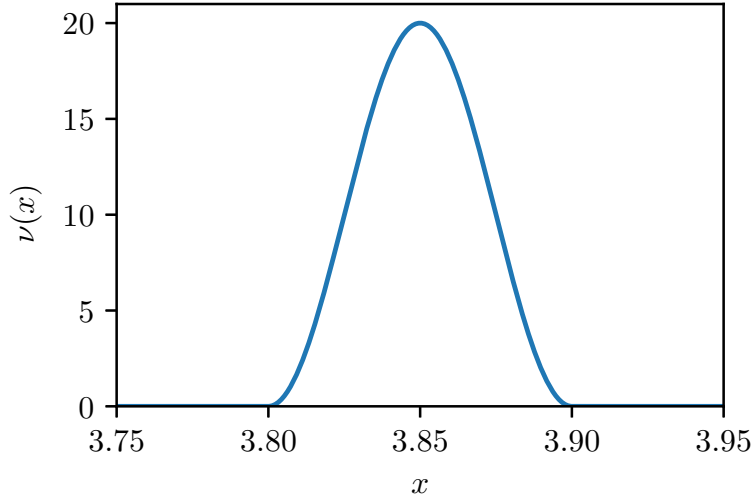}
	\caption{Probability density $\nu(x)$ of the raised cosine distribution supported on $[3.8,3.9]$.}
	\label{f.raisedcosine}
\end{figure}

\change{
We study perturbations of the form 	
\[ a^{(j)} = a_0^{(j)} +  \epsilon a_1^{(j)}, \]
where 
the $a_0^{(j)}$ are sampled from a $C^1$ compactly supported distribution $\nu(a_0) da_0$ and the $a_1^{(j)}$ are sampled from a compactly supported distribution $\nu(a_1|a_0) da_1$ depending smoothly on $a_0$. As we argue below, smoothness of $\nu$ is crucial for the existence of linear and higher-order response. For concreteness, we choose $\nu(a_0) da_0$ to be the raised cosine distribution supported on the interval $[3.8,3.9]$:  
\[ \nu(x) = \frac{\mathbf{1}_{[3.8,3.9]}(x) }{0.2} \left(1 + \cos \frac{x-3.85}{0.05}\pi \right), \] which is depicted in Figure \ref{f.raisedcosine}. The raised cosine distribution is bell-shaped and has a second derivative with bounded variation. This degree of smoothness implies that  our model exhibits cubic (third-order) response, as will be shown in the next section. Furthermore, for the numerical simulations in Section~\ref{s.homo} and \ref{s.det} we choose $a_1^{(j)} = 1$ for all $j$.
}

\change{
We shall consider two cases, $\gamma=\tfrac{1}{2}$ and $\gamma=1$ corresponding to a diffusive scaling limit and a deterministic scaling limit, respectively. In the thermodynamic limit $M\to \infty$, we will see that in the former case the microscopic driving term $Z_n$ converges to a stochastic process $\zeta_n$ in the macroscopic dynamics, whereas in the latter case $Z_n$ converges to a constant. This, in conjunction with heterogeneously distributed microscopic parameters $a^{(j)}$, leads asymptotically to macroscopic linear response in the former case, and a failure of linear response in the latter.
}

\section{$\gamma=\frac{1}{2}$: Weak coupling with diffusive limit}
\label{s.homo}
We now justify LRT for the high-dimensional system (\ref{e.Qn})-(\ref{e.qn})  with $\gamma=\tfrac{1}{2}$. This is done in two steps. We first show that the dynamics of the macroscopic variable $Q$ is diffusive. The invariant measure of this diffusive process depends on the integrated effect of the microscopic variables for a specific configuration of the parameters $a^{(j)}$. In a second step we establish conditions on the parameter distribution $\nu(a)$ for the logistic map parameters of the microscopic sub-system which allow for expectation values of an observable of the resolved state to vary smoothly with the perturbation size $\eps$.

\begin{figure}
\centering
\includegraphics[width=0.7\linewidth]{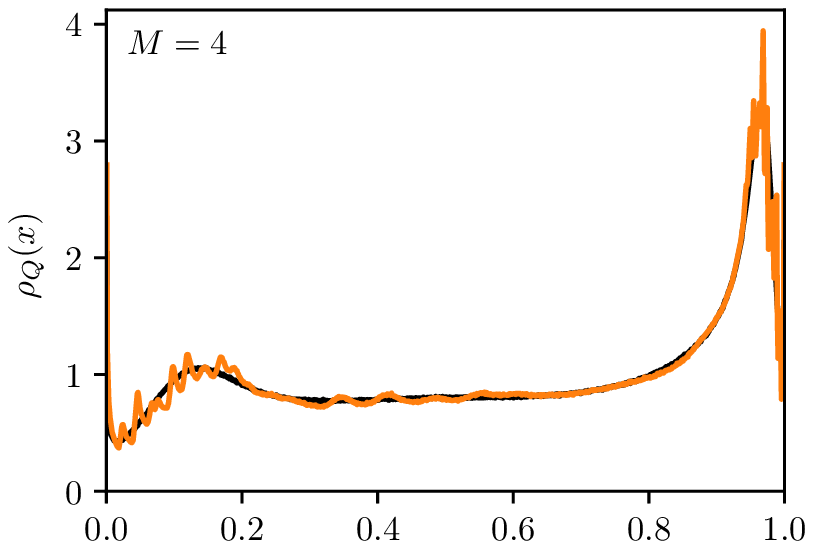}
\includegraphics[width=0.7\linewidth]{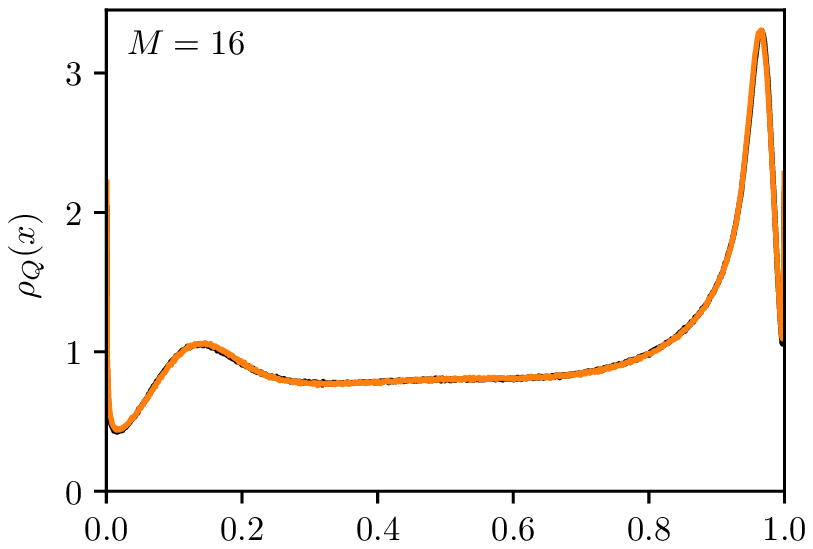}
\caption{Empirical probability density $\rho_Q(x)$ (orange) of the macroscopic variable $Q$ for $\gamma=\tfrac{1}{2}$ as estimated from simulations of the original deterministic system (\ref{e.Qn})-(\ref{e.qn}) for different values of the size $M$ of the microscopic sub-system. Top: $M=4$. Bottom: $M=16$. The continuous black line depicts the invariant density of the stochastic limit system (\ref{e.Qnstoch}). We used $A_0=3.91$, $A_1=0.05$ and $\eps=0$.}
\label{f.g0p5_pdf}
\end{figure}

\begin{figure}
\centering
\includegraphics[width=0.7\linewidth]{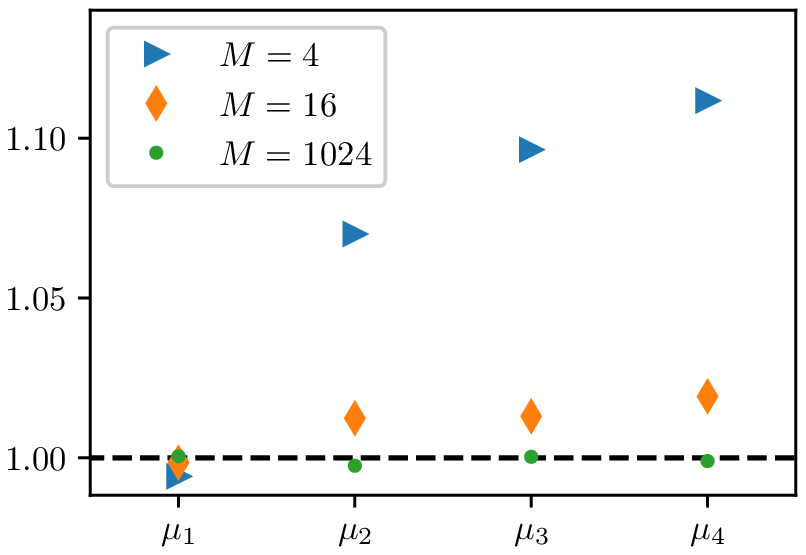}
\caption{First four centred moments $\mu_{i}$, $i=1,\cdots,4$, of the macroscopic variable $Q$ for $\gamma=\tfrac{1}{2}$ as estimated from simulations of the original deterministic system (\ref{e.Qn})-(\ref{e.qn}) for fixed time $n=6$ for several values of the size  of the microscopic sub-system: $M=4$ (blue triangles), $M=16$ (orange diamonds) and $M=1024$ (green dots). We depict the moments scaled by the respective moments of the stochastic limit system (\ref{e.Qnstoch}) so that the asymptotic limit is $1$ for all moments. Parameters as in Fig.~\ref{f.g0p5_pdf}.}
\label{f.g0p5_mom}
\end{figure}

We begin by considering the unperturbed case $\eps=0$ and show that the macroscopic variable $Q$ asymptotically satisfies a stochastic limit system in the thermodynamic limit $M\to\infty$ when $\gamma=\tfrac{1}{2}$. We consider driving terms $Z_n$ with mean-zero functions $\phi(\cdot;a)$, \change{$\E[\phi(\cdot;a)]=0$, where the average is with respect to the invariant measure of the unresolved microscopic variable for fixed parameter $a=a_0^{(j)}+\eps a_1^{(j)}$. (Whenever we consider averages for fixed parameters rather than for fixed $\eps$ we drop the superscript of $\E$).} The driving term $Z_n$ contains a sum over independent identically distributed random variables for each time $n$. Hence, for $\gamma=\tfrac{1}{2}$, the central limit theorem assures that the driving term $Z_n$ converges to a random Gaussian variable $\zeta_n\sim{\mathcal{N}}(0,\sigma^2)$ with $\sigma^2=\langle \phi^2\rangle$, where the angular brackets denote the average over the measure of the logistic map parameters $\nu(da)$. 
Moreover, in discrete time the $\zeta_n$ define a stationary Gaussian stochastic process -- a moving average process of infinite order --  which is (subject to continuity assumptions \cite{KarlinTaylor75}) uniquely defined by its mean and its covariance $R(m)$. The covariance is readily determined as
\begin{align}
R(m)={\mathrm{cov}}(\zeta_n,\zeta_{n+m})
=\lim_{M\to\infty}\frac{1}{M}\sum_{j=1}^M \E[\phi^{(j)}_0\phi^{(j)}_m]  
 = \langle \E[\phi_0\phi_m] \rangle .
\label{e.covzetan}
\end{align}
The process $Q_n$ hence converges weakly to the stochastic process defined by
\begin{align}
\Q_{n+1} = (A_0+A_1\zeta_n)\,\Q_n(1-\Q_n).
\label{e.Qnstoch}
\end{align}
Figure~\ref{f.g0p5_pdf} illustrates the convergence of the deterministic map (\ref{e.Qn})-(\ref{e.qn}) to the stochastic limit system (\ref{e.Qnstoch}) in distribution by comparing the respective empirical measures for several values $M$ of the size of the microscopic sub-system. The microscopic dynamics is run unperturbed with $\eps=0$. Here we chose the mean-zero (conditional on the parameter $a$) functions $\phi(x,a)=x^2- \left( a x (1-x)\right)^2 $ to generate the driving sum $Z_n$. We used a time series of $N=4\times 10^7$ and determined the empirical measure of the full system (\ref{e.Qn})-(\ref{e.qn}) by binning using $1000$ bins. Details on how to determine the statistics of the limiting diffusive system (\ref{e.Qnstoch}) are given in Appendix \ref{appb}. It is remarkable that with only $M=16$ microscopic variables the eye can barely distinguish the empirical density from the density of the diffusive limit equation (\ref{e.Qnstoch}). We further show convergence of the first four moments of $Q$ when increasing $M$ in Figure~\ref{f.g0p5_mom}. It is seen that for accurate convergence of higher order moments to the values of their stochastic limiting equation (\ref{e.Qnstoch}) larger system sizes $M$ are required. 

After having established that the dynamics of the macroscopic variable $Q$ is diffusive, we now establish in a second step that the associated invariant measure and expectation values of macroscopic observables depend smoothly on $\eps$. It is pertinent to stress that the mere existence of a stochastic limit does not imply LRT. \change{We remark that the existence of the stochastic limit is in line with the Gallavotti-Cohen hypothesis; however, this is insufficient for LRT.} Consider, for example, the case when each unresolved microscopic variable $q^{(j)}$ evolves according to the logistic map with the same parameters $a^{(j)}\equiv {\rm{const}}$, differing only in the initial conditions drawn from the invariant measure. The limit system would still be a stochastic system due to the randomness in the initial conditions, but LRT would not be valid when homogeneously perturbing the unresolved scales. Crucial for the validity of LRT is that the parameters $a_0^{(j)}, a_1^{(j)}$ are identically independently distributed ({\em{i.i.d.}}) random variables, sampled from a distribution $\nu(da_0,da_1)$ with a regularity property that we now derive. 

For microscale variables with parameter $a$, consider the expectation of an observable $\Phi_a=\E[\phi_0]$ or $\Phi_a=\E[\phi_0\phi_m]$, and consider its average over the microscopic dynamics $\langle \Phi\rangle_{\eps}=\int \Phi_{a_0+\eps a_1} \nu(a_0,a_1)da_0da_1
$. Changing variables $\alpha=a_0+\eps a_1$ we find $\langle \Phi\rangle_{\eps}=\int \Phi_\alpha \nu(\alpha-\eps a_1,a_1)\, d\alpha da_1 $, and hence 
\begin{align*}
\frac{d}{d\eps} \langle \Phi\rangle_\eps = -\int a_1 \Phi_{a_0+\eps a_1}\,\frac{\partial}{\partial a_0}\nu(a_0,a_1)\, da_0 da_1 . 
\end{align*}
This is well-defined provided that $a_1 \frac{\partial}{\partial a_0} \nu(da_0,da_1)$ is integrable: if so, the statistics of $\zeta_n$ vary smoothly with respect to $\eps$. A particular case of this is when $a_0$ and $a_1$ are independently distributed and the marginal density of $a_0$ is of bounded variation. It is readily seen that to achieve higher-order response, say of order $\ell$, (weak) derivatives of order $\ell$ must be defined. This can be achieved if $a_0$ and $a_1$ are drawn independently from a distribution $\nu$ with a marginal distribution $\nu(a_0)$ in Sobolev space $W^{\ell,1}$.
  
\begin{figure*}[htbp]
\centering
\includegraphics[width=1\linewidth]{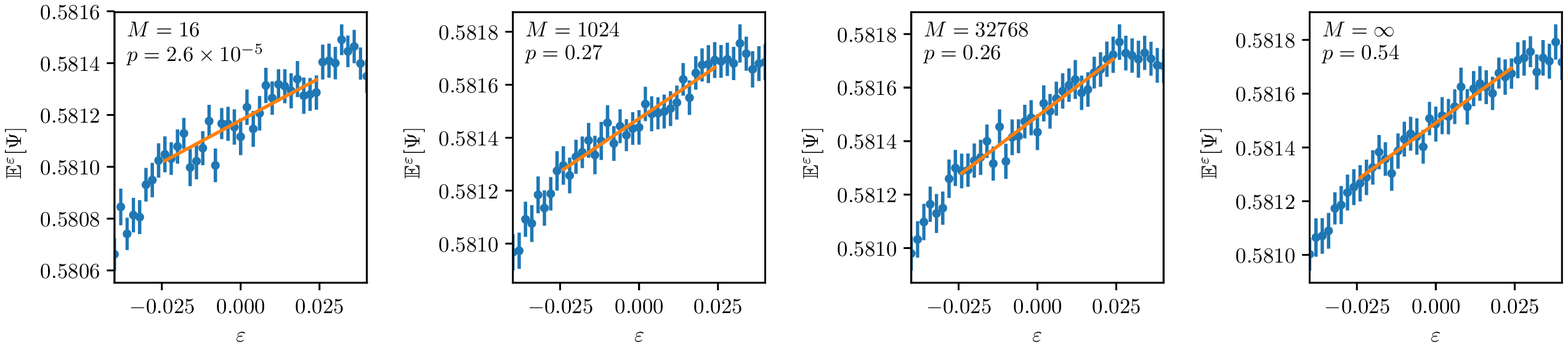}
\caption{Linear response of an observable $\Psi(Q)=Q$ for the deterministic system (\ref{e.Qn})-(\ref{e.qn}) 
for $\gamma=\tfrac{1}{2}$ for different values of the size $M$ of the microscopic sub-system. (a): $M=16$. (b): $M=1024$. (c) $M=32768$. (d): Stochastic limit system (\ref{e.Qnstoch}). All experiments used a time series of length $N=2\times10^5$. The error bars were estimated from $K=200$ realizations differing in the initial conditions. We used $A_0=3.91$, $A_1=0.05$.}
\label{f.g0p5_LRT}
\end{figure*}

We present in Figure~\ref{f.g0p5_LRT} results of the linear response for an observable $\Psi(Q)=Q$. The microscopic sub-system is perturbed homogeneously with $a_1^{(j)}=1$ for all $j$. It is clearly seen that the perturbation $\eps$ induces a smooth change in the observable for large $M$, indicative of the validity of LRT. We employ here the test for linear response developed in \cite{GottwaldEtAl16} and report the $p$-values testing the null hypothesis of linear response. We compute averages for several values of $\eps$ from long simulations of length $N=5\times 10^6$. The error bars shown in Figure~\ref{f.g0p5_LRT} are estimated from $K=200$ realizations differing in the initial conditions of the microscopic variables. For completeness we provide a description and justification of the test in Appendix \ref{appa}. For small values of $M=16$ the $p$-value is $\O(10^{-5})$, rejecting the null hypothesis of linear response, whereas for $M=2^{10}$ the $p$-value is $0.27$, consistent with linear response. We also show results of the linear response for the stochastic limit system (\ref{e.Qnstoch}), illustrating that the thermodynamic limit implies linear response with a $p$-value of $p=0.54$. Note that although the invariant density of the resolved degree of freedom $Q$ has sufficiently converged to the invariant density of the stochastic limit system (\ref{e.Qnstoch}) for $M=16$ (cf. Figure~\ref{f.g0p5_pdf}), this size is not sufficiently large to assure linear response. 

In Figure \ref{f.app1} in Appendix \ref{appa} we present results for cubic response for the same simulations which gave rise to Figure~\ref{f.g0p5_LRT}. Cubic response is valid for the stochastic limiting system (\ref{e.Qnstoch}) because the raised cosine distribution, which was chosen for the distribution $\nu$ in the simulation, lies in $W^{3,1}$.

\section{$\gamma=1$: Weak coupling with deterministic limit}
\label{s.det}
\begin{figure}
\centering
\includegraphics[width=0.7\linewidth]{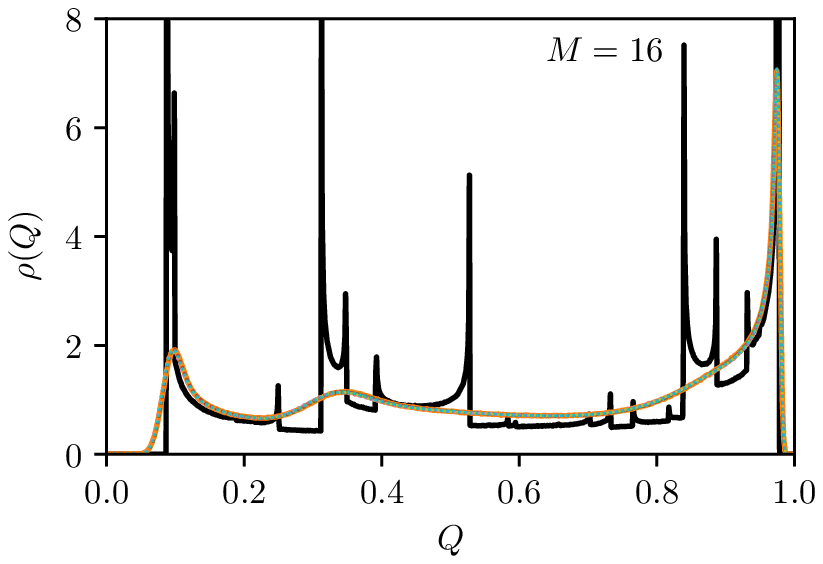}
\includegraphics[width=0.7\linewidth]{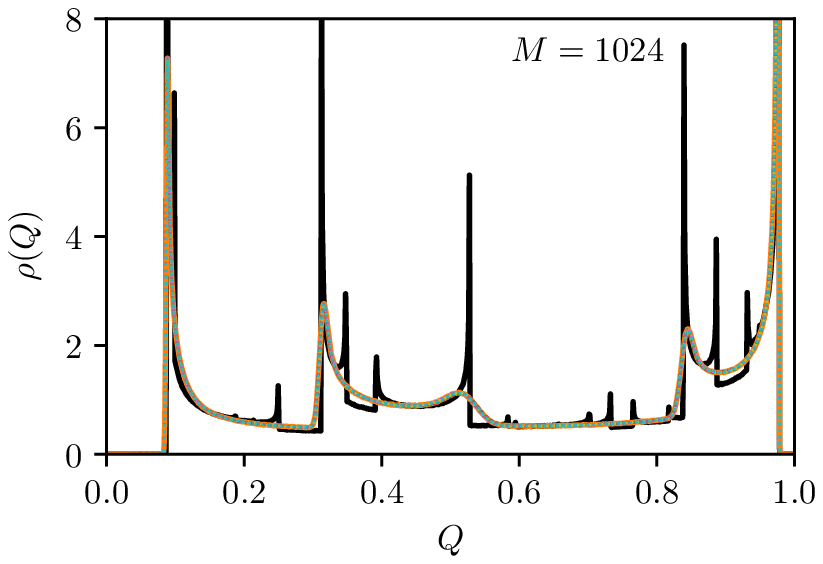}
\caption{Empirical probability density $\rho_Q(x)$ (orange line) of the macroscopic variable $Q$ for $\gamma=1$ as estimated from simulations of the original deterministic system (\ref{e.Qn})-(\ref{e.qn}) for different values of the size $M$ of the microscopic sub-system. Top: $M=16$. Bottom: $M=1024$. The continuous black line depicts the invariant density of the deterministic logistic map limit system (\ref{e.Qnstoch}); the thin dotted lines, which are indistinguishable from $\rho_Q(x)$, 
represent invariant densities of the logistic map (\ref{e.Qn}) with the stochastic driving $Z_n$ given by (\ref{e.Zndet}) for various realisations of $\eta$.   
%
We used $A_0=3.847$, $A_1=0.147$ and $\eps=0$.}
\label{f.g1_pdf}
\end{figure}

We begin again by considering the unperturbed case $\eps=0$. In the case $\gamma=1$ we consider the driving term $Z_n$ generated by a function $\phi$ with non-vanishing mean and consider $\phi(x,a)=x^2$.  Since each unresolved degree of freedom generates an invariant measure, for $\gamma=1$ the driving variable $Z_n$ converges to a constant according to the law of large numbers with $Z_n \to C=\langle \E[\phi] \rangle$. 
In the thermodynamic limit therefore the limiting equation is a deterministic logistic map
\begin{align}
\Q_{n+1} = A\,\Q_n(1-\Q_n)
\label{e.Qndet}
\end{align}
with $A = A_0 + C A_1$. Figure~\ref{f.g1_pdf} illustrates the convergence of the invariant measure of the deterministic map (\ref{e.Qn})-(\ref{e.qn}) to the averaged deterministic limit system (\ref{e.Qnstoch}) in distribution upon increasing the size $M$ of the microscopic sub-system. We used again a time series of $N=4\times 10^7$ and determined the empirical measure by binning using $1000$ bins. We see that for $M=1024$ convergence to the rough limiting invariant measure of the deterministic logistic map (\ref{e.Qndet}) with its narrow peaks has not been fully achieved. This is due to finite sample size $M$. \change{There are two averages requiring the limit $M\to\infty$: the average $\E[\cdot]$ with respect to the invariant density of the microscopic logistic dynamics and the average $\langle \cdot\rangle$ with respect to the distribution of the parameters $a^{(j)}$ of the modified logistic dynamics. Each of those averages is associated with their own finite size correction which are captured by the central limit theorem. Up to $\mathcal{O}(1/\sqrt{M})$ we have
\begin{align}
Z_n &= \frac{1}{M}\sum_{j=1}^{M}\phi^{(j)}(q_n^{(j)})= \frac{1}{M}\sum_{j=1}^{M}\E[\phi^{(j)}] + \frac{1}{\sqrt{M}}\zeta_n
\nonumber \\
&=\langle \E[\phi] \rangle + \frac{1}{\sqrt{M}}\eta + \frac{1}{\sqrt{M}}\zeta_n.
\label{e.Zndet}
\end{align}
Here $\zeta_n$ is again the mean-zero Gaussian process with covariance matrix $R(m)$ defined in (\ref{e.covzetan}), and for fixed $\eps$, $\eta$ is a Gaussian variable with $\eta\sim{\mathcal{N}}(0,\langle\E^\eps[\phi]^2\rangle - \langle \E^\eps[\phi]\rangle^2)$. In the context where $\eps$ varies, $\eta$ in (\ref{e.Zndet}) can be understood as a random function of $\eps$, having a mean-zero Gaussian distribution with covariance 
\[\langle \eta^\eps \eta^{\eps'} \rangle = \langle \E^\eps[\phi] \E^{\eps'}[\phi] \rangle - \langle \E^\eps[\phi] \rangle \langle \E^{\eps'}[\phi] \rangle.\]
In general, $\eta$ is non-differentiable which implies that LRT is violated for macroscopic observables $\Psi(Q)$, even for the random finite-size driver $Z_n$ given by (\ref{e.Zndet}). However, if the variation in $\E[\phi]$ over the parameter values sampled by $\nu$ is small by comparison with the typical variance $R(0)=\E[(\phi-\E[\phi])^2]$ for these parameters (e.g. if the support of $\nu$ is sufficiently small), then the small, rough contribution of $\frac{1}{\sqrt{M}}\eta$ to the response of $\Psi(Q)$ is dominated by the (linear) response generated by $\langle \E[\phi] \rangle + \frac{1}{\sqrt{M}}\zeta_n$. We remark, however, that if the support of $\nu$ is too small and the parameters are therefore less heterogeneous, LRT is only valid for a small range of perturbation sizes $\eps$. 
}

\change{To illustrate the role of finite size effects, we present in Figure~\ref{f.g1_pdf} also results of simulations of the logistic map (\ref{e.Qn}) with $Z_n$ stochastically generated by (\ref{e.Zndet}), 
mimicking random finite size effects in approximating the deterministic limit $Z_n = \langle \E [\phi] \rangle$. It is seen that for finite $M$ the peaks are smoothed by sampling noise, and the random logistic map reproduces the invariant density of the macroscopic variable $Q$ of the full deterministic model driven by the microscopic dynamics.}\\

\begin{figure*}[htbp]
\centering
\includegraphics[width=1\linewidth]{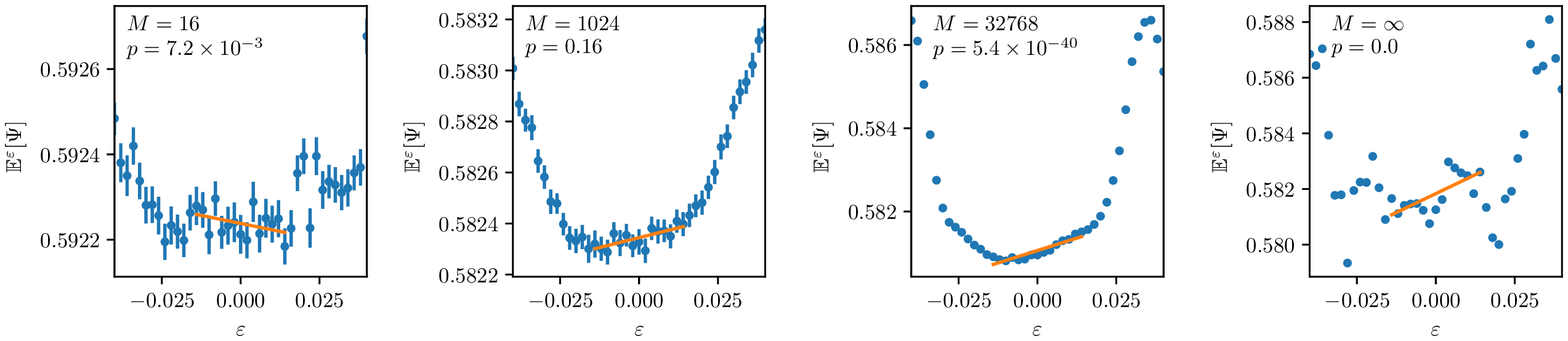}
\caption{Linear response of an observable $\Psi(Q)=Q$ for the deterministic system (\ref{e.Qn})-(\ref{e.qn}) for $\gamma=1$ for different values of the size $M$ of the microscopic sub-system. (a): $M=16$. (b): $M=1024$. (c) $M=32768$. (d): Deterministic limit system (\ref{e.Qndet}). All experiments used a time series of length $N=2\times10^5$. The error bars were estimated from $K=200$ realizations differing in the initial conditions. We used $A_0=3.847$ and $A_1=0.147$.}
\label{f.g1_LRT}
\end{figure*}

Given that the thermodynamic limit system is deterministic, one might be tempted to conclude that linear response is not valid. Figure~\ref{f.g1_LRT} shows the linear response as a function of perturbation $\eps$ for several values of the microscopic sub-system size $M$. For small values of $M$ LRT is clearly violated with a $p$-value of $\O(10^{-3})$, as expected. For very large values of $M=2^{15}$ LRT is violated with a $p$-value of $\O(10^{-40})$, consistent with the LRT-violating deterministic limit system (\ref{e.Qndet}). Remarkably and maybe surprisingly, decreasing the size $M$ from these very large values to intermediate values of $M=1024$ we observe that the numerical results are consistent with LRT and the $p$-value increases dramatically to around $0.16$. \change{This can be explained by the finite size corrections (\ref{e.Zndet}) to the deterministic limit $Z_n = \langle \E[\phi] \rangle$ provided by the central limit theorem. We note that the p-value for $M=1024$ indicates marginal evidence in favour of breakdown of LRT associated with the (small) contribution of the non-differentiable $\eta$ term.} 
Just as in the $\gamma = \tfrac{1}{2}$ case it is necessary for LRT to hold in the case of finite sample size, that the parameters $a^{(j)}$ are inhomogeneously distributed with a sufficiently smooth distribution $\nu(a)$.\\

\begin{figure}
\centering
\includegraphics[width=0.7\linewidth]{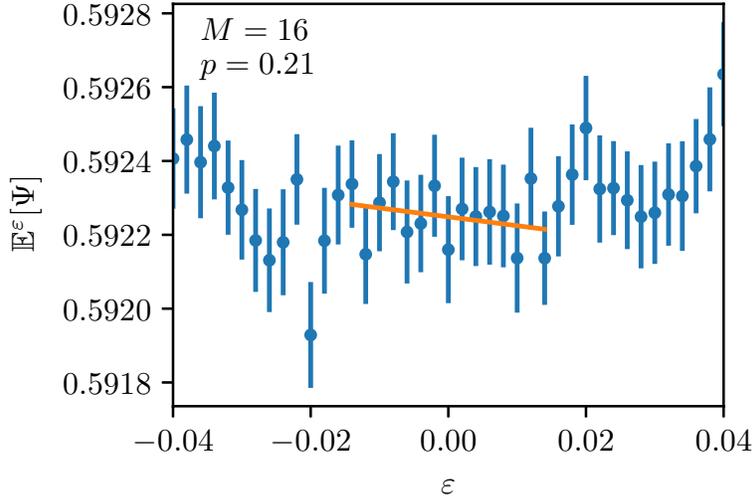}
\caption{Linear response of an observable $\Psi(Q)=Q$ for the deterministic system (\ref{e.Qn})-(\ref{e.qn}) for $\gamma=1$ with $M=16$ estimated from a time series of length $N=2\times 10^4$. The error bars were estimated from $K=200$ realizations differing in the initial conditions. We used $A_0=3.847$ and $A_1=0.147$.}
\label{f.g1_LRT_finitesize}
\end{figure}

In \cite{GottwaldEtAl16} it was found that even if a system does not obey linear response one might not be able to reject the null hypothesis of linear response with sufficient statistical significance when the data length $N$ of the time series is not sufficiently long. In Figure~\ref{f.g1_LRT_finitesize} we show the linear response as a function of $\eps$ for a microscopic sub-system of size $M=16$ for $N=2\times 10^4$. While for $N=2\times 10^5$ linear response was rejected with $p=7.2\times 10^{-3}$, linear response is now consistent with the given data with a $p$-value of $p=0.21$. It is found that decreasing the length of the time series allows for a larger range in the perturbation size $\eps$ for which linear response is consistent with the data. 

\section{Discussion and outlook}
\label{s.discussion}
We have shown that macroscopic observables in high-dimensional deterministic dynamical systems which consist of unresolved microscopic variables weakly coupled to macroscopic resolved variables may obey linear response theory even if each of the microscopic units individually violate LRT. We showed that in the case when the thermodynamic limit of an infinitely large microscopic sub-system leads to a stochastic limit equation for the macroscopic resolved variables, linear response theory can be justified for macroscopic observables. In case when the thermodynamic limit is deterministic we showed that for a finite microscopic sub-system, the limiting dynamics has a stochastic correction which again allows for linear response. We established that the existence of a stochastic limit system is not sufficient to assure LRT, and an additional assumption on the distribution of the parameters of the microscopic sub-system is needed in the case when the microscopic variables are not respecting linear response. In this case, we require the parameters of the microscopic sub-system to be heterogeneous with a smooth distribution of their parameters. The degree of the smoothness directly determines the polynomial order of the response. For example, if the parameters of the unresolved degrees of freedom $q^{(j)}$ were chosen to be all equal and the initial conditions $q_0^{(j)}$ were chosen from the invariant measure, the macroscopic variable $Q$ still obeyed a stochastic limit for $\gamma=\tfrac{1}{2}$, but LRT would clearly be violated upon a homogeneous perturbation of the microscopic sub-system. If the microscopic variables obey linear response, for example with uniformly expanding maps, this condition on the parameter distribution is not necessary.

We considered here the worst case scenario where the dynamics of both the macroscopic and the unresolved degrees of freedom on their own violate LRT. In the numerical simulations we studied the effect of perturbing the parameters of the unresolved microscopic variables. \change{We remark that if perturbations of the macroscopic variable $Q$ were considered with $A=A_0+\varepsilon \delta A$, LRT would be valid for $\gamma=\tfrac{1}{2}$ since the limiting system is stochastic \cite{HairerMajda10} (and also for $\gamma=1$ when finite size effects are non-negligible).} Rather than considering a macroscopic variable weakly coupled to a micrscopic sub-system consisting of non-conservative logistic maps, one may instead consider the case of a traditional heat bath consisting of an infinite collection of harmonic oscillators with randomly chosen initial conditions which are weakly coupled to a distinguished resolved degree of freedom. The limiting stochastic properties of the associated $Z_n$ was established rigorously in \cite{FordEtAl65,Zwanzig73,FordEtAl87,StuartWarren99,KupfermanEtAl02,GivonEtAl04} using trigonometric approximation of Gaussian noise \cite{Kahane85}. In this case, if weakly coupled to the macroscopic variable $Q$ which evolves according to the logistic dynamics (\ref{e.Qn}), we would obtain similar results as for the case considered here. 

We have treated here the case of weakly coupled systems. It is well known that stochastic limit systems also occur in multi-scale dynamics where the central limit theorem is realized by summing up many fast chaotic degrees of freedom in one slow time unit \cite{GivonEtAl04,Dolgopyat04,MelbourneStuart11,GottwaldMelbourne13c,KellyMelbourne14,DeSimoiLiverani15,DeSimoiLiverani16}. We expect analogous results in this case. As in the case of weak coupling considered here, the heterogeneity in the parameter distribution of the fast system is essential.


\section*{Acknowledgements}
GAG is partially supported by ARC, grant DP180101385. CW is supported by an Australian Government Research Training Program (RTP) Scholarship.

\appendix

\section{Model reduction for chaotic microscopic sub-systems}
\label{appb}

This appendix describes how to compute the statistics of the stochastic limiting system Eqn (\ref{e.Qnstoch}) for $\gamma=\tfrac{1}{2}$, which we recall here
\begin{align}
\Q_{n+1} = (A_0+A_1\zeta_n)\,\Q_n(1-\Q_n),
\label{e.Qnstoch2}
\end{align}
for the deterministic limiting system Eqn (\ref{e.Qndet}) for $\gamma=1$, which we recall here
\begin{align}
\Q_{n+1} = A\,\Q_n(1-\Q_n)
\label{e.Qndet2}
\end{align}
with $A = A_0 + \langle \E[\phi] \rangle A_1$, and for the stochastic finite-size system 
\begin{align}
\Q_{n+1} = A\,\Q_n(1-\Q_n)
\label{e.Qndet3}
\end{align}
with $A = A_0 +Z_nA_1$ where $Z_n$ is given by Eqn (\ref{e.Zndet}), which is recalled here as
\begin{align}
Z_n  = \langle \E[\phi] \rangle + \frac{1}{\sqrt{M}}\eta + \frac{1}{\sqrt{M}}\zeta_n.
\label{e.Zndet2}
\end{align}
The random variable $\eta$ accounts for the random variation in the selection of the parameters $a^{(j)}$ and the random process $\zeta_n$ accounts for the dynamics of the microscopic variables. (However, as can be seen from Figure 
\ref{f.g1_pdf}, setting $\eta \equiv 0$, i.e. replacing it with its expectation, gives a remarkably good approximation of the invariant measure, at least in the system we consider.)

In order to simulate these systems we need to estimate $\langle \E[\phi] \rangle$ and, for the stochastic systems, also  $R(m) = \langle \E[\phi_0 \phi_m] \rangle - (\E [\phi])^2, m \in \mathbb{N}$. We describe first how we estimate these parameters from Monte Carlo simulations of the logistic map, and then describe how we sample the stochastic process $\zeta_n$ with the covariance parameters given by $R(m)$.

\subsection{Estimating parameters}
We need to estimate the expectation values for $K$ perturbation sizes $\eps_i$ with $i=1,\cdots,K$. Since we set here $a^{(j)}=1$ for all microscopic variables, at each $\eps_i$ we write the averages over the microscopic dynamics as
\begin{equation} 
\langle \E^{\varepsilon_i} [\phi ]\rangle = \int_{\R} \E^{\alpha} [\phi(\cdot,\alpha)] \, \nu(\alpha-\varepsilon_i) d\alpha 
\label{e.Logexp}
\end{equation}
and
\begin{equation} 
\langle \E [\phi_0 \phi_m] \rangle_{\varepsilon_i} = \int_{\R} \E^{\alpha} [\phi(x_0,\alpha) \phi(x_m,\alpha)m]  \, \nu(\alpha-\varepsilon_i) d\alpha \label{e.Logcov}
\end{equation}
for $i = 1, \ldots, K$ and  $m = 1, \ldots, \infty$, where $\nu$ is the density function of the logistic map parameters and is chosen here as the raised cosine distribution
\[ 
\nu(a) = \mathbf{1}_{[3.8,3.9]} \frac{1}{0.1} \left(1+\cos\left(\frac{a-3.85}{0.05}\pi\right)\right).
\]
From now on it is understood that all observables, expectations and so on are for a fixed parameter $\alpha$: we therefore drop the $\alpha$ and $(j)$ superscripts for ease of exposition.

We use a trapezoidal rule to estimate the integrals in (\ref{e.Logexp}) and (\ref{e.Logcov}), using a grid of $30,001$ values of the logistic map parameters $\alpha$ evenly spaced on $[3.7,4.0]$ (to allow for the support of $\nu$ as well as the range of the perturbation). This is used for each $\varepsilon_i$.

The expectations (\ref{e.Logexp}) and (\ref{e.Logcov}) can be entirely determined by simulations of a standard logistic map without coupling to the expanding $r$-dynamics: Denote by $\varphi_n = \phi(x_n,a_0+\eps_i)$ such that $x_{n+1}= (a_0+\eps_i) x_n(1-x_n)$ with $x_0 = q_0$. The logistic dynamics of the $q$ are augmented by $r$-dynamics so that at any time step the $q$ will with equal probability either advance according to the logistic map or remain unchanged. The invariant measure of $q$ is therefore identical to the one supported by a logistic map with the same parameter $\alpha$; hence $\E [\phi] = \E [\varphi]$.\\ To estimate the averages of the auto-correlations  (\ref{e.Logcov}) we define $N(m)$ as the number of evolution steps of the $q$-dynamics up to physical time $m$ which were done according to the logistic map (i.e. discarding all those instances when the $r$-dynamics forces $q$ not to vary). Note that $N(m)$ has a binomial distribution $N(m)\sim {\mathrm{B}}(m,\tfrac{1}{2})$. Hence by definition we have
\[ 
\phi(q_m) = \varphi_{N(m)},
\]
and we can write
\begin{align*}
\E [\phi_0 \phi_m] &= \E [\varphi_0 \varphi_{N(m)}] \\
&=  \sum_{i=0}^m 2^{-i} {m \choose i} \E [\varphi_0 \varphi_i]. 
\end{align*}
For regular values of $\alpha$, when the logistic map $x_n$ with parameter $\alpha$ has a stable periodic orbit, calculating the stable periodic orbit allows for an accurate evaluation of the expectation. We use the database of periodic windows given in \cite{Galias17} to identify regular points and stable periodic orbits.

For chaotic values of $\alpha$ we estimate expectations and lag-correlations of the logistic map with parameter $\alpha$ via Monte-Carlo simulation of the logistic map $x_{n}$, using $10$ separate initialisations with $399168$ time steps each. This number of time steps was chosen as it has a large number of prime factors, and therefore will give more accurate estimates for short periodic windows outside the database, or for chaotic values where the acim has multiple connected components (i.e., $f$ is not mixing but $f^p$ is for some $p>1$).

\subsection{Sampling the stochastic process $\zeta_n$}

The limiting process $\zeta_n$ is a stationary Gaussian process given by lag-covariance function $R(m)$. Assuming sufficiently fast decay of the lag-covariance function, we can write this process as a moving-average process of infinite order
\[ \zeta_n = \sum_{m=0}^\infty \beta_m X_{n-m} \]
with a deterministic sequence  $(\beta_m)_{m\in \mathbb{N}}\in \ell_2$ and {\em i.i.d.} standard normal random variables $X_n$.

The moving average coefficients $\beta_m$ and the covariance function $R_m$ are related by
\[ 
R(m) = \sum_{k=0}^\infty \beta_k \beta_{m+k}. 
\]
The coefficients can be extracted from the covariance function via the generating functions 
\[
\mathcal{B}(z) := \sum_{m=0}^\infty \beta_m z^m 
\]
and 
\[
\mathcal{R}(z) := \sum_{m=-\infty}^\infty R(|m|) z^m, 
\]
for which the relation $\mathcal{R}(z) = \mathcal{B}(z) \mathcal{B}(z^{-1})$ holds. If we restrict to the complex unit circle, setting $z = e^{i\theta}$, we find that $ \mathcal{R}(e^{i\theta}) = \mathcal{B}(e^{i\theta})\mathcal{B}(e^{-i\theta}) = |\mathcal{B}(e^{i\theta})|^2$ since $\beta_m\in \mathbb{R}$. Assuming that $\mathcal{R}(e^{i\theta}) \neq 0$, we have that
\[ 
\frac{1}{2} \log \mathcal{R}(e^{i\theta}) = \Re \log \mathcal{B}(e^{i\theta}),
\]
and hence, we can write, using that the $\beta_m$ are real, 

\[
\log \mathcal{B}(z) = \sum_{m=0}^\infty b_m z^m
\]
with $b_m \in \R$. The $b_m$ may be calculated via Fourier cosine transform using that
\[ 
\frac{1}{2} \log \mathcal{R}(e^{i\theta}) = \sum_{m=0}^\infty b_m \cos m\theta.
\]
The $b_m$ coefficients allow one to evaluate $\mathcal{B}(e^{i\theta})$, from which the moving average coefficients $\beta_m$ are obtained via an additional Fourier transform.


\section{Testing for linear response in finite time series}\label{appa}

We summarize here briefly the quantitative goodness-of-fit test for the detectability of linear response introduced in \cite{GottwaldEtAl16}. The test quantifies the statistical significance of an observed linear response in time series of finite size. 

Given a family of chaotic maps $f_\varepsilon$ that may or may not obey linear response, we test for linear response at some reference state with parameter $\varepsilon=\varepsilon_0$. Hence we seek to examine the linear dependency of the response of a bounded and continuous observable
\begin{align}
\delta \Psi = \E^\eps[ \Psi] - \E^{\eps_0} [\Psi]
\end{align}
in terms of the perturbation size $\eps$. To test for linearity we consider $K > 2$ different values of the perturbation parameter $\varepsilon_1,\dotsc,\varepsilon_K$, and sample $N$ consecutive values from the perturbed maps yielding the time series $x^i_n=f\ueps{i}(x^i_{n-1})$ for each $i=1,\dotsc, K$ and $n=1,\dotsc, N$. The initial conditions $x^i_0$ are distributed according to the physical measure associated with $f\ueps{i}$. The lengths of the time series $N$ is chosen that for each $i=1,\dotsc,K$ the corresponding autocorrelation function has sufficiently decayed, i.e.  we choose $N_i \gg \tau_{\Psi,\varepsilon_i}$, where $\tau_{\Psi,\varepsilon_i}$ is the $1/e$-folding time of the autocorrelation function of $\Psi$ under the dynamics $f_{\varepsilon_i}$. For simplicity, we choose $N_i=N$ for all $i$ in the following.\\

It is well known that for a large class of chaotic dynamical systems, the sample averages of the observations
\begin{align}
\bar{\Psi}_i=\frac{1}{N}\sum_{n=1}^{N}\Psi(x^i_n)
\label{e.Aavgs}
\end{align}
obey the central limit theorem and are distributed asymptotically as $\mathcal{N}\left(\pav{\Psi}{{i}}, \sigma_i^2/N\right)$ \cite{melbournelectures,ColletEckmann07} with 
\begin{align}
\bar{\Psi}_i \approx \E^{\eps_i}[\Psi]  + \frac{\sigma_i}{\sqrt{N}} \xi_i \; ,
\label{e.clt}
\end{align}
for $i = 1,\dotsc,K$ and {\em{iid}} noise $\xi_i \sim \mathcal{N}(0,I)$. The variances $\sigma_i^2$ are given by the Green-Kubo formula as an infinite sum of lag-correlations of $f\ueps{i}$ as
\begin{align} 
\sigma_i^2 = C_0(\Psi,\Psi) + 2 \sum_{j=1}^{\infty} C_j(\Psi,\Psi),
\label{e.covar}
\end{align}
where the correlation function $C_n$ between two observables $\Psi$ and $\Omega$ is defined as
\begin{align*}
C_n(\Psi,\Omega) = \paveps{{\Psi \hphantom{;}\Omega\circ f^n}} -\, \paveps{\Psi} \paveps{\Omega}.
\end{align*}
The variances can be efficiently estimated numerically as a Monte-Carlo estimate from (\ref{e.clt}) for large $N$. The results in this work were obtained with $N=40\times 10^6$.

In the case the dynamical system has linear response and provided the perturbations $\delta \varepsilon_i=\varepsilon_i-\varepsilon_0$ are sufficiently small, the following statistical model holds for $\bar{\Psi}_i$ (with $o(\delta\varepsilon_i)$ error)
\begin{align}
\bar{\Psi}_i  = \alpha_0 + \alpha_1 \, \delta\varepsilon_i + \frac{\sigma_i}{\sqrt{N}} \xi_i \, ,
\label{e.linmod}
\end{align}
with $\alpha_0=\E^{\eps_0}[\Psi]$ and $\alpha_1= \E^{\eps_0}[\Psi]^\prime$ for the unperturbed reference state with $\varepsilon=\varepsilon_0$. It is pertinent to mention that the $\xi_i$ are independent since the samples from each perturbed system are generated independently. 

The parameters $\alpha_0$ and $\alpha_1$ of the model (\ref{e.linmod}) can be determined from time series by means of a weighted least squares fit and we obtain
\begin{align*}
\left(\begin{array}{c}
\hat \alpha_0\\
\hat \alpha_1
\end{array}\right) = (D^T D)^{- 1} D^T Y
\end{align*}
with the design matrix
\begin{align*}
D = \left(\begin{array}{cc}
1 / \sigma_1 & \delta \varepsilon_1 / \sigma_1\\
\vdots & \vdots\\
1 / \sigma_K & \delta \varepsilon_K / \sigma_K
\end{array}\right)\, ,
\end{align*}
and the vector of scaled observations
\begin{align*}
Y = \left(\begin{array}{c}
\bar{\Psi}_1 / \sigma_1\\
\vdots\\
\bar{\Psi}_K / \sigma_K
\end{array}\right)\, .
\end{align*}

Testing for validity of linear response then amounts to testing whether the actual observations could have been generated from the linear model (\ref{e.linmod}) with normally distributed errors $\xi_i\sim {\mathcal{N}}(0,I)$. To do so we choose a Pearson $\chi^2$-test to test the goodness-of-fit with statistics
\begin{align}
\chi^2 
&= N\, \sum_{i=1}^K \left(Y_i - \frac{1}{\sigma_i}\left(\hat \alpha_0+\hat \alpha_1\varepsilon_i\right)\right)^2 
\nonumber
\\
&= N\, Y^T(I-H)Y,
\label{e.chisq}
\end{align}
where the idempotent matrix 
\begin{align*}
H = D (D^T D)^{-1} D^T
\end{align*} 
maps scaled observations $Y$ to their linear fits, i.e. $HY = D(\hat \alpha_0\;\;\hat \alpha_1)^T$ \cite{BoxHunterHunter}.\\
If the response of the underlying dynamical system is linear, $\chi^2$ has a $\chi^2$-distribution with $K-2$ degrees of freedom and expectation value $\E \chi^2_{K-2} = K-2$. 
Hence a measure for the breakdown of linear response can be quantified as the difference between the $\chi^2$ test statistic for the scaled observations $Y_i= \bar{\Psi}_i / \sigma_i$ and the expectation of the test statistic under the null hypothesis of linear response
\begin{align} 
\q= \frac{1}{N}\left( \chi^2 -\E \chi^2_{K-2} \right).
\label{e.Echi2_0}
\end{align}
The central limit theorem (\ref{e.clt}) holds independent of the existence of linear response and can be used to obtain expressions for the mean and variance of the breakdown parameter. Defining $W$ as the vector with components $W_i=\E^{\eps_i}[\Psi]/\sigma_i$, the mean is calculated as
\begin{align} 
\E \q&= \frac{1}{N}\left(\E\chi^2 -\E \chi^2_{K-2}\right)
\nonumber
\\ 
&=  \E\bigg((W + \frac{1}{\sqrt{N}}\bxi)^T(I-H)(W + \frac{1}{\sqrt{N}}\bxi) - \frac{1}{N}\E \chi^2_{K-2}\bigg)
\nonumber
\\
&= \|W-HW\|^2,
\label{e.Echi2}
\end{align}
where we used that $H$ is idempotent. Hence $\q$ is a random variable whose expected value measures the difference between the actual response $\E^{\eps_i}[\Psi]$ and an assumed linear response $\alpha_0+\alpha_1\varepsilon_i$ as calculated via least square regression. The mean of the breakdown parameter $\E \q$ is non-negative and is zero only for $W = HW$, i.e. if the observations stem from a dynamical system obeying linear response. The variance of the breakdown parameter $\q$ is calculated as
\begin{align*} 
\V \q &= \E \bigg((W + \frac{1}{\sqrt{N}}\bxi)^T(I-H)(W + \frac{1}{\sqrt{N}}\bxi) -\frac{K-2}{N} - \E \q \bigg)^2\\
&= \frac{1}{N} \E\left( \bxi^T (I-H) (2W + \frac{1}{\sqrt{N}}\bxi)  - \frac{K-2}{\sqrt{N}}  \right)^2.
\end{align*}
Note that $\V\q\to 0$ for $N\to \infty$, and hence $\q$ is a consistent estimator for the mismatch $\E\q = \|W-HW\|^2$. In numerical experiments it is practical to consider Monte-Carlo estimates of the mismatch over realizations $\q_j$ differing in their initial condition and set
\begin{align}
\hat\q=\frac{1}{K} \sum_{j=1}^K\q_j\, .
\end{align}

To make statements about the statistical significance of whether an observed time series of length $N$ is classified as obeying linear response or not, we introduce a $p$-value testing the null hypothesis of linear response. Let us consider the case when a dynamical system does not obey linear response, i.e. when $\E\q\neq 0$. Using Chebyshev's inequality we have that for all $b < N\E\q$,
\begin{align*}
P(N\q < b) &\leq P(|\q-\E\q| > \E\q - b/N)\\
&\leq \frac{\V(\q)}{(\E\q - b/N)^2}\; .
\end{align*}
Since, as we have shown above, $\V \q\to 0$ as $N\to\infty$, we conclude that $N \q \to \infty$ in probability as $N\to\infty$. Hence, if $F$ is the cumulative distribution function of the $\chi^2_{K-2}$ distribution, the $p$-value obtained using the $\chi^2$- test,
\begin{align}
p  = 1 - F(\chi^2) = 1- F(K-2 +N\q),
\label{e.p}
\end{align}
converges quickly in probability to zero as $N\to\infty$ \cite{BoxHunterHunter}. This implies that the probability of falsely accepting the null hypothesis of linear response at any significance level can be made arbitrarily small for sufficiently large data length $N$.\\

For completeness (although not used in this work) we show that one can define a threshold value $\q_\alpha$ for the observed random variable $\hat\q$ such that if $\hat \q >\q_\alpha$ the null hypothesis of linear response is rejected with significance level $\alpha$ (i.e. with probability $1-\alpha$); given a specified significance level $\alpha$ the threshold value can be defined as
\begin{align}
\label{e.qalpha}
\q_\alpha = \frac{1}{N} \left(F^{-1}(1-\alpha)-(K-2)\right)\, .
\end{align}
It is clear from this that the detectability of breakdown of linear response crucially depends on the amount of available data. As $N\to \infty$, a breakdown will always become detectable at any specified significance level $\alpha$. Conversely, if the mismatch $\E\q$ between the true response of the dynamical system and the linear response is too small and there is an insufficient amount of data available, the actual response will be swamped by the sampling noise, and one will not be able to detect the breakdown of linear response with a reasonable significance level.\\

\begin{figure*}
	\centering
	\includegraphics[width=1\linewidth]{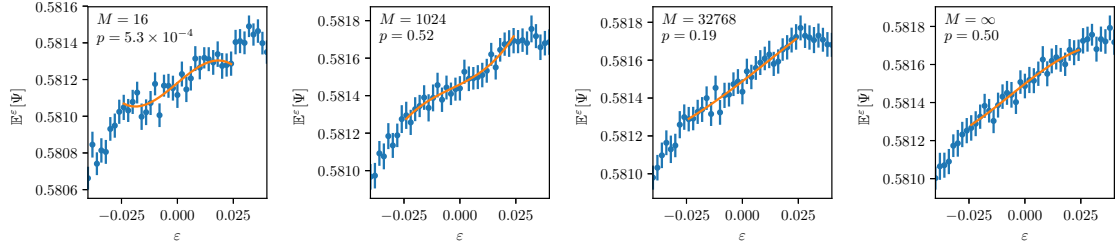}
	\caption{Cubic response of an observable $\Psi(Q)=Q$ for the deterministic system (1)--(3) for $\gamma=\tfrac{1}{2}$. (a): $M=16$. (b): $M=1024$. (c) $M=32768$. (d): Stochastic limit system (5). All experiments used a time series of length $N=2\times10^5$. The error bars are estimated from $K=200$ realizations differing in the initial conditions. We used $A_0=3.91$, $A_1=0.05$.}
	\label{f.app1}
\end{figure*}

It is possible to extend this test to probe higher order response. To test for $\ell$th order response we add terms $\sum_{j=2}^\ell \alpha_j \delta\eps_i^j$ to our statistical model (\ref{e.linmod}) and then employ higher-order regression (i.e. augmenting the design matrix $D$). In Figure~\ref{f.app1} we show results for the same numerical simulations as in Figure~\ref{f.g0p5_LRT}, but now showing cubic response rather than linear response. We recall that we can expect cubic response due to the distribution density of the logistic map parameters of the microscopic sub-system being three times continuously differentiable. As for linear response, the null-hypothesis of cubic response can be rejected for small values of the system size $M$ but cannot be rejected for sufficiently large values of $M$. For $M=16$ the test yields a $p$-value of $2.1\times 10^{-5}$. For $M=1024$ the $p$-value is $0.19$ consistent with cubic response.\\

For more details on the test the interested reader is referred to \cite{GottwaldEtAl16}.\\

\bibliographystyle{siam}
\bibliography{bibliography}

\end{document}